\documentclass[MSNbibl,number,citesort,seceqn,dvips]{arxbj}
\usepackage{upgreek}
\usepackage{mathbh}

\aid{0}
\volume{18}
\issue{3}
\pubyear{2012}
\firstpage{783}
\lastpage{802}
\doi{10.3150/11-BEJ367}

\makeatletter
\renewcommand{\emptyset}{\varnothing}
\renewcommand{\epsilon}{\varepsilon}
\renewcommand{\pi}{\uppi}
\def\Rset{\mathbb R}
\def\bbr{\mathbb R}
\def\Zset{\mathbb Z}
\def\esp{\mathbb E}
\def\pr{\mathbb P}
\def\cov{\operatorname{cov}}
\def\rme{\mathrm{e}}
\def\rmi{\mathrm{i}}
\def\rmd{\mathrm{d}}
\def\1{\mathbh{1}}

\def\mce{\mathcal{E}}
\def\mcf{\mathcal{F}}
\def\mcw{\mathcal{W}}
\def\mcx{{\mathcal X}}
\def\mcz{{\mathcal Z}}

\newtheorem{theorem}{Theorem}
\newtheorem{lemma}[theorem]{Lemma}
\newtheorem{corollary}[theorem]{Corollary}
\newremark{remark}{Remark}
\newproclaim{assumption}{Assumption}
\newcommand{\eqref}[1]{(\ref{#1})}
\makeatother

\begin{document}
\begin{frontmatter}

\vspace*{-3pt}
\title{Function-indexed empirical processes
based on an infinite source Poisson transmission stream}
\runtitle{The empirical process}

\begin{aug}
\author[1]{\fnms{Fran\c{c}ois} \snm{Roueff}\corref{}\thanksref{1}\ead[label=e1]{roueff@telecom-paristech.fr}},
\author[2]{\fnms{Gennady} \snm{Samorodnitsky}\thanksref{2}\ead[label=e2]{gennady@orie.cornell.edu}} \and
\author[3]{\fnms{Philippe} \snm{Soulier}\thanksref{3}\ead[label=e3]{philippe.soulier@u-paris10.fr}}
\runauthor{F. Roueff, G. Samorodnitsky and P. Soulier}
\address[1]{Institut TELECOM, Telecom Paristech, CNRS LTCI,
46, rue Barrault,
75634 Paris Cedex 13, France. \printead{e1}}
\address[2]{School of Operations Research and Information Engineering,
and Department of Statistical Science,
Cornell University,
Ithaca, NY 14853, USA. \printead{e2}}
\address[3]{Modal'X, Universit{\'e} Paris Ouest Nanterre,
200 avenue de la R{\'e}publique,
92000 Nanterre, France.\printead{e3}}
\end{aug}

\received{\smonth{10} \syear{2010}}
\revised{\smonth{2} \syear{2011}}

%
\begin{abstract}
We study the asymptotic behavior of empirical processes generated by
measurable bounded functions of an infinite source Poisson transmission
process when the session length have infinite variance. In spite of the
boundedness of the function, the normalized fluctuations of such an empirical
process converge to a non-Gaussian stable process. This phenomenon can be
viewed as caused by the long-range dependence in the transmission process.
Completing previous results on the empirical mean of similar types of
processes, our results on nonlinear bounded functions exhibit the influence
of the limit transmission rate distribution at high session lengths on the
asymptotic behavior of the empirical process. As an illustration, we apply
the main result to estimation of the distribution function of the steady
state value of the transmission process.
\end{abstract}

%

\end{frontmatter}
%

\section{Introduction}

We consider the infinite source Poisson transmission process defined by
%
\begin{equation} \label{eqdefgginfini}
X(t) = \sum_{\ell\in\Zset} W_\ell  \1_{\{ \Gamma_\ell\leq t
< \Gamma_\ell+ Y_\ell\}},\qquad t\in\Rset  ,
\end{equation}
where the triples $\{(\Gamma_\ell,Y_\ell,W_\ell), \ell\in\Zset\}
$ of
session arrival times, durations and transmission rates satisfy the
following assumption.
\begin{assumption}\label{hypoisp}
\renewcommand\theenumi{(\roman{enumi})}
\renewcommand\labelenumi{\theenumi}
\begin{enumerate}[(iii)]
\item\label{itempoisson} The arrival times $\{\Gamma_\ell, \ell
\in\Zset\}$ are the
points of a homogeneous Poisson process on the real line with
intensity $\lambda$, indexed in such a way that $\cdots<
\Gamma_{-2} < \Gamma_{-1} < \Gamma_{0} <0<\Gamma_{1}< \Gamma_2 <
\cdots.$\vadjust{\goodbreak}
\item\label{itemmarques} The durations and transmission rates $\{(Y,W),
(Y_\ell,W_\ell), \ell\in\Zset\}$ are independent and identically
distributed random pairs with values in
$(0,\infty)\times[0,\infty)$ and independent of the arrival times
$\{\Gamma_\ell, \ell\in\Zset\}$. The random variables $W_j$ are
positive with a positive probability. The session lengths $Y_j$ have
finite expectation and infinite variance.
\item\label{itemjointe} There exist a measure $\nu$ on
$(0,\infty]\times[0,\infty]$ such that
$\nu((1,\infty]\times[0,\infty]) = 1$ and, as $n\to\infty$,
\[
n\pr \biggl( \biggl(\frac{Y}{a(n)},W \biggr) \in\cdot \biggr)
\stackrel v \longrightarrow\nu ,
\]
where $\stackrel v \rightarrow$ denotes vague convergence on
$(0,\infty]\times[0,\infty]$, and $a$ is the left continuous
inverse $(1/\bar F)^\leftarrow$ of $1/\bar F$. Here $F$ is the
distribution function of $Y$, and $\bar F = 1-F$ is the corresponding
survival function. The relatively compact sets of $(0,\infty]\times
[0,\infty]$
are all sets contained in $[\epsilon,\infty]\times[0,\infty]$ for
some positive
$\epsilon$, see Resnick \cite{resnick2007}, Chapter~3.
\end{enumerate}
\end{assumption}

Assumption~\ref{hypoisp}\ref{itemjointe} implies several things,
listed below. See Heffernan and Resnick \cite{heffernanresnick2007}.
\begin{itemize}
\item The
survival function $\bar F$ is regularly varying with index
$-\alpha$ for some $\alpha>0$. The function $a$
is then regularly varying with index $1/\alpha$.
\item The limiting measure $\nu$ is a product measure:
%
\begin{equation} \label{eqdef-G}
\nu= \nu_\alpha\times G ,
\end{equation}
where $\nu_\alpha$ is a measure on $(0,\infty)$ satisfying
$\nu_\alpha((x,\infty)) = x^{-\alpha}$ for all $x>0$, and $G$ is a
probability measure on $[0,\infty]$.
\item We have the following weak convergence on $[0,\infty]$, as $t\to
\infty$,
%
\begin{equation} \label{eqconv-to-G}
\pr (W\in\cdot| Y>t ) \stackrel w \longrightarrow G .
\end{equation}
\end{itemize}

We will assume that the exponent $\alpha$ satisfies
%
\begin{equation} \label{ealpha}
1<\alpha<2 .
\end{equation}

Under Assumption~\ref{hypoisp}, the process (\ref{eqdefgginfini}) is
well defined and stationary, see, e.g., Fay, Roueff and Soulier \cite{fayroueffsoulier2007}.
Under additional moment assumptions, it is shown in this reference
that the autocovariance function of the process
$X$ is regularly varying at infinity with index $2H-2\in(-1,0)$,
where $H=(3-\alpha)/2$. Such slow rate of decay of the covariance
function is often associated with long range dependence.

We are interested in studying the large time behavior of the empirical
process
%
\begin{equation} \label{eqmcz}
\mathcal J_T(\phi) = \int_0^T \phi ( X_h(s) ) \, \rmd s, \qquad
 T>0,
\end{equation}
where $h>0$, $X_h(s) = \{X(s+t)   , 0 \leq t \leq h\}$, and $\phi$ is
a real valued measurable function defined on the space
$\mathcal{D}([0,h])$ endowed with the $J_1$ topology, see, for instance,
Kallenberg~\cite{kallenberg2002}. We notice that the $\mathcal{D}([0,h])$-valued
stochastic process $ ( X_h(s),   s\in[0,T] )$ is continuous
in probability and, hence, has a measurable version, see
Cohn \cite{cohn1972}. In particular, $\mathcal J_T(\phi)$ above is a
well defined random variable, as long as the function $\phi$ satisfies
appropriate integrability assumptions, for example, when the function
$\phi$
is bounded.

The case $h=0$ and $\phi(x) = x$ has been considered in
Mikosch \textit{et al}. \cite{mikoschresnickrootzenstegeman2002} with $W_i\equiv1$ and by
Maulik, Resnick and Rootz{\'e}n \cite{maulikresnickrootzen2002} in the present context of possible
dependence between the session lengths and the rewards (transmission
rates). These references consider the case where the intensity of the point
process of arrivals is possibly increasing, which gives rise to the slow
growth/fast growth dichotomy. In the slow growth case, which includes
the case
of constant intensity, the limit of the partial sum process is a L\'evy stable
process, whereas in the fast growth case, the limiting process is the
fractional Brownian motion with Hurst index $H = (3-\alpha)/2$. Here, we
consider a fixed intensity for the sessions arrival rate, hence are restricted
to the slow growth case. On the other hand, we take~$\phi$ arbitrary (but
bounded) and thus obtain what appears to be the first result on the asymptotic
behavior of the empirical process for this type of long range dependent shot
noise process. The limit process depends on the intensity $\lambda$,
the tail
exponent $\alpha$ and the limit transmission rate distribution $G$ defined
in (\ref{eqconv-to-G}). As an illustration, we apply the main result
to the
estimation of the distribution function of the steady state value of the
transmission process. Moreover, we allow $h>0$. Other potential
applications of
our main result (e.g., to estimation of the multivariate distribution function)
can be handled in a~similar way, but we do not pursue them in this paper.

Our main result is stated as a functional central limit theorem in the Skorohod
$M_1$ topology. A convergence result in this topology was obtained in
Resnick and van~den Berg~\cite{resnickvandenberg2000a} for a similar traffic model, but with
$h=0$ and
$\phi(x) = x$. Our result can be viewed as a heavy traffic
approximation of
the content of a fluid queue fed with input $\phi(X(s))$. It shows, in
particular, that even for $\phi$ bounded (e.g., with $\phi(x)=x\wedge
b$ with
$b$ denoting a maximal allowed bandwidth), the fluctuations of the asymptotic
approximation of the queue content have an infinite variance. See also
Resnick and van~den Berg \cite{resnickvandenberg2000a}, Section~5.

\section{Notation and preliminary results}
\label{secnotat-prel-results}
We now introduce some notation and derive certain useful properties of
the empirical process (\ref{eqmcz}) stated in several lemmas whose
proofs are provided in Section~\ref{secproof-lemmas}.

We employ the usual queuing terminology: a time point $t$ is said to
belong to a busy period if $X(t)>0$; it belongs to an idle period
otherwise. A cycle consists of a busy period and the subsequent idle
period.

The following facts about M/G/$\infty$ queues will be useful, see
Hall \cite{hall1988}. Under Assumption~\ref{hypoisp}\ref{itempoisson}
and \ref{itemmarques}, one can define the sequence $\{S_j, j\in
\Zset\}$ of
the successive starting times of the cycles such that
$\cdots <S_{-2}<S_{-1}<0<S_0<S_1<\cdots $. Define the cycle lengths
$C_j=S_j-S_{j-1}$ for all $j\in\Zset$. Hence, $S_0$ is the starting
time of
the first complete cycle starting after time $0$ (note that $S_0$ may
or may
not be equal to the first Poisson arrival after time 0), and $S_n =
S_0+\sum_{j=1}^{n} C_j$. The cycle form a regenerative sequence in the sense
that $\{(C_j, X(\cdot+S_{j-1}) \1_{[0,C_j)}), j \geq1 \}$ is an i.i.d.
sequence of random pairs with values in $(0,\infty)\times\mathcal
D([0,\infty))$. Moreover, we have
%
\begin{equation}
\label{eqcycle-tail}
\esp[C_1] = \rme^{\lambda\esp[Y]}/\lambda .
\end{equation}
The following result provides the tail behavior of $C_1$. It is proved in
Section~\ref{secproof-lemmas}.
\begin{lemma} \label{lemtailC}
Suppose that Assumption~\ref{hypoisp} holds. Then
$C_1$ has a regularly varying tail with index $\alpha$ and
%
\begin{equation}
\lim_{t\to\infty} t  \pr\bigl(C_1>a(t)x\bigr) = \rme^{\lambda\esp[Y]}
x^{-\alpha}   . \label{ecyclel}
\end{equation}
\end{lemma}

Let $\phi$ be a measurable function defined on $\mathcal{D}([0,h])$,
satisfying
appropriate integrability conditions for the integral in \eqref
{eqmcz} to be
well defined (e.g., bounded). We decompose~$\mathcal J_T(\phi)$ using
the cycles
defined above. Let us denote
%
\begin{equation}\label{eqZj}
Z_j(\phi) = \int_{S_{j-1}}^{S_{j}} \phi(X_h(s)) \, \rmd s,   \qquad
j=1,2,\ldots.
\end{equation}
Then $(Z_j(\phi))_{j\geq1}$ is a stationary sequence, but, if $h>0$,
it is not
an i.i.d. sequence. Nevertheless, it is easy to see that it is strongly
mixing. Define the sigma-fields $\mcf_j = \sigma(Z_k(\phi), 1\leq k
\leq j)$ and
$\mathcal G_j = \sigma(Z_k(\phi), k > j)$ and mixing coefficients
$(\alpha_k)_{k\geq1}$ by
\[
\alpha_k = 2 \sup\{|\cov(\1_A,\1_B)|,  A \in\mcf_j, B \in
\mathcal G_{j+k}  ,  j\geq1\} .
\]
Let $j,k\geq1$, $A\in\mcf_j$ and $B\in\mathcal G_{j+k}$. Denote $U=
\1_A-\pr(A)$ and $V = \1_B-\pr(B)$. Then
\begin{eqnarray*}
|\cov(\1_A,\1_B)| \leq\pr(S_{j+k}-S_j \leq h) + \bigl|\esp\bigl[ UV
\1_{\{S_{j+k}-S_{j}>h\}}\bigr]\bigr|.
\end{eqnarray*}
Observe that $U\1_{\{S_{j+k}-S_{j}>h\}}$ is $\sigma\{X(S_{j+k}-t),
t>0\}$-measurable, $V$ is $\sigma\{X(S_{j+k}+t),
t\geq0\}$-measurable and that by the regenerative property, these two
sigma-fields are independent. Thus, $\esp[ UV
\1_{\{S_{j+k}-S_{j}>h\}}]=0$ and we obtain, for all $k\geq1$,
%
\begin{equation}
\label{eqmixing}
\alpha_k  \leq2 \sup_{j\geq1} \pr(S_{j+k}-S_{j}\leq h) \leq
2\sup_{j\geq1} \pr\bigl(\max(C_{j+1},\ldots,C_{j+k}) \leq h\bigr) =
2F_C(h)^{k}   ,
\end{equation}
where $F_C$ denotes the distribution function of $C_1$. Since
$F_C(h)<1$ for any $h$, the mixing coefficients $\alpha_k$ decay
exponentially fast, independently of $\phi$. This property will be a~key ingredient to the proof of our result since it implies that, in
many aspects, the
sequence~$Z_j(\phi)$ has the same asymptotic properties as an i.i.d.
sequence.

Let $\mce(\cdot,\phi)$ be the function defined on $[0,\infty)$ by
%
\begin{equation}\label{eqPhiDef}
\mce(w,\phi) = \esp\bigl[\phi\bigl(w+X_h(0)\bigr)\bigr] ,
\end{equation}
whenever the latter expectation is well defined, which is always the
case if
$\phi$ is bounded. In that case, by Fubini's theorem,
$\mce(\cdot,\phi)$ is a measurable function.
It follows from the elementary renewal theorem that
$\esp[Z_j(\phi)] = \esp[\phi(X_h(0))]\esp[C_1]$. This identity is stated
formally in the following lemma, which also contains another result
that will
be needed later.

\begin{lemma}
\label{lemtruly-iid}
Suppose that Assumption~\ref{hypoisp} holds. Let $h\geq0$ and
$\phi$ be a bounded measurable function defined on
$\mathcal{D}([0,h])$. We have
%
\begin{equation}
\label{eqZmean}
\esp [Z_1(\phi) ]=\mce(0,\phi)\esp[C_1] =
\esp[\phi(X_h(0))]\esp[C_1]  .
\end{equation}
Moreover, for any $p\in(1,\alpha)$, there exists a constant $C>0$ and
a positive function $g$ depending neither on $\phi$ nor on $T$ such
that $g(x)\to0$ as $x\to\infty$ and
%
\begin{equation} \label{equnifBound}
\pr \Bigl(\sup_{t\in[0,T]} |\mathcal J_{t}(\phi) - \esp
[\mathcal
J_{t}(\phi)] |>x
\|\phi\|_\infty \Bigr) \leq C T^{1-p} + C T  x^{-p}+ g(x)   .\vspace*{-2pt}
\end{equation}
\end{lemma}

For all $\epsilon,t>0$, let $N_{\epsilon,t}$ be the number of
sessions of
length greater than $\epsilon a(t)$ arriving and ending within the first
complete cycle $[S_0,S_1)$. Further, we let $Y_{\epsilon,t}$ be the
length of
the first such session starting at or after $S_0$ with length greater than
$\epsilon a(t)$ and let $\Gamma_{\epsilon,t}$ and~$W_{\epsilon,t}$ be,
correspondingly, its starting time and the transmission rate. The following
lemma shows that, when $N_{\epsilon,t}\geq1$, the process $\{\phi
(X_h(s)) ,
s\in[S_0,S_1)\}$ can be, in certain sense, approximated by the step function
$\{\mce(W_{\epsilon,t},\phi)\1_{[\Gamma_{\epsilon,t},\Gamma
_{\epsilon,t}+Y_{\epsilon,t})}(s)
 ,   s \in[S_0,S_1)\}$. (Note that by definition, if $N_{\epsilon
,t}\geq1$,
then $S_0 \leq\Gamma_{\epsilon,t}< \Gamma_{\epsilon,t}+Y_{\epsilon
,t} \leq
S_1$.)\vspace*{-2pt}

\begin{lemma}\label{lemM1}
Suppose that Assumption~\ref{hypoisp} holds. Let $h\geq0$ and
$\phi$ be a bounded measurable function defined on
$\mathcal{D}([0,h])$. Let $\eta>0$. We have, for all $\epsilon>0$
sufficiently
small,
%
\begin{equation} \label{eqM1cycle}
\pr\biggl ( \sup_{v\in[S_0,S_1]}  \biggl| \int_{S_{0}}^{v}
\bigl\{\phi(X_h(s)) - \mce(W_{\epsilon,t}, \phi)
\1_{[\Gamma_{\epsilon,t},\Gamma_{\epsilon,t} +
Y_{\epsilon,t})}(s)\bigr\} \, \rmd s  \biggr| > \eta a(t)  ;
N_{\epsilon,t} \geq1  \biggr) = \mathrm{o}(t^{-1})  .\vspace*{-2pt}
\end{equation}
\end{lemma}

Let $\mcw$ be a closed subset of $[0,\infty]$ such that
$\pr(W\in\mcw)=1$. (Note that by (\ref{eqconv-to-G}) this implies
$G(\mcw)=1$.) We introduce the following assumption.\vspace*{-2pt}
\begin{assumption}\label{hypocty}
We have
%
\begin{equation}
\label{eqdiscpointshavenulmeasure}
G (D(\mce(\cdot,\phi),\mcw) )=0 ,
\end{equation}
where $D(\mce(\cdot,\phi),\mcw)$ denotes the set of discontinuity
points of the function $\mce(\cdot,\phi)$
restricted to $\mcw\cap[0,\infty)$, and containing the point $\infty$
if $\infty\in\mcw$ and $\mce(w,\phi)$ does not converge as
$w\to\infty$ with $w\in\mcw$. (The notation $\mce(\infty,\phi)$, when
used in the sequel, refers to the continuous extension of
$\mce(w,\phi)$, and will be used only when such an extension
exists.)\looseness=-1\vspace*{-2pt}
\end{assumption}

\begin{remark}\label{remdiscCase}
If the distribution of $W$ is supported on a closed set consisting
of isolated points in $[0,\infty)$ (which would be the case, for
instance, if $W$ was a nonnegative integer-valued random variable),
then $D(\mce(\cdot,\phi),\mcw)$ is either empty or equal to
$\{\infty\}$. In the latter case, if $G(\{\infty\})=0$, then
Assumption~\ref{hypocty} is verified.\vspace*{-2pt}
\end{remark}

The next lemma, which may be of independent interest, states the
multivariate regular variation property of the empirical process over
a cycle.\vspace*{-2pt}

\begin{lemma} \label{lemjoint-tails}
Suppose that Assumption~\ref{hypoisp} holds. Let $h\geq0$ and
$\phi_1,\ldots,\phi_d$ be bounded measurable functions
defined on $\mathcal{D}([0,h])$ satisfying Assumption\vadjust{\goodbreak}
\ref{hypocty}
with $G$ defined by (\ref{eqdef-G}). With $\mce(w,\phi_i) =
\esp[\phi_i(w+X_h(0))]$, $i=1,\ldots, d$, $w\geq0$, we let
\begin{eqnarray*}
\mathbf Z = \biggl [\int_{S_0}^{S_1} \phi_1(X_h(s)) \, \rmd s,\ldots,
\int_{S_0}^{S_1} \phi_d(X_h(s)) \, \rmd s \biggr]^{\mathrm{T}}  .
\end{eqnarray*}
Then $\mathbf Z$ is multivariate regularly varying with index
$\alpha$. More precisely, the following vague convergence holds on
$[-\infty,\infty]^d\setminus\{0\}$ as $t\to\infty$,
%
\begin{equation}
\label{eqqueueZmultdim}
t\pr \biggl( \frac{\mathbf Z}{a(t)} \in\cdot \biggr) \stackrel v
\longrightarrow
\rme^{\lambda\esp[Y]}   \int_{y=0}^\infty\pr\bigl (y[\mce
(W^*,\phi_1),\ldots,
\mce(W^*,\phi_d)]^{\mathrm{T}}\in\cdot \bigr) \alpha y^{-\alpha-1}\,\rmd y
  ,
\end{equation}
where $W^*$ is a random variable with values in $[0,\infty]$ and
distribution $G$.
\end{lemma}

\section{Main result}

As observed in Resnick and van~den Berg \cite{resnickvandenberg2000a}, since the limit is
discontinuous, the convergence of the sequence of processes
$\{\mathcal Z_{T}(\phi,t), t\geq0\}$ in
Theorem~\ref{theostable-constant} cannot hold in $\mathcal
D([0,\infty))$ endowed
with the topology induced by Skorohod's $J_1$
distance. We shall prove that the convergence holds in $\mathcal
D([0,\infty))$
endowed with the topology induced by Skorohod's $M_1$ distance.

\begin{theorem}
\label{theostable-constant}
Suppose that Assumption~\ref{hypoisp} holds. Let $h\geq0$ and
$\phi$ be a bounded measurable function on $\mathcal{D}([0,h])$
satisfying Assumption \ref{hypocty} with $G$ defined
by (\ref{eqdef-G}). Then, as $T\to\infty$, the sequence of
processes $\mathcal Z_T(\phi,\cdot)$ defined by
%
\begin{equation} \label{eqemp-process}
\mathcal Z_T(\phi,u) = \frac1{a(T)} \int_0^{Tu} \{ \phi(X_h(s)) -
\esp[\phi(X_h(0))] \} \, \rmd s , \qquad   u\geq0 ,
\end{equation}
converges weakly in $\mathcal D([0,\infty))$ endowed with the $M_1$
topology to a strictly $\alpha$-stable L\'evy motion $ (
\Lambda(\phi,u),   u\geq0 )$ satisfying
%
\begin{equation} \label{eqlimit-process}
\esp \mathrm{e}^{\mathrm{i}t\Lambda(\phi,u)} = \exp \bigl\{ -u|t|^\alpha\lambda
c_{\alpha}
\esp | \mce(W^*,\phi) - \mce(0,\phi) |^\alpha\{1 - \rmi
\beta\operatorname{sgn}(t) \tan(\pi\alpha/2)\}  \bigr\}
\end{equation}
for $u\geq0$ and $t\in\bbr$, where $ c_\alpha =
-\Gamma(1-\alpha)\cos(\pi\alpha/2)$, $W^*$ is
as in Lemma \ref{lemjoint-tails}, and
\[
\beta= \frac{\esp [ | \mce(W^*,\phi) - \mce(0,\phi
) |^\alpha
\operatorname{sgn} ( \mce(W^*,\phi) -
\mce(0,\phi) ) ]}{\esp | \mce(W^*,\phi) -
\mce(0,\phi) |^\alpha}  .
\]
\end{theorem}

\begin{remark}
For applications of Theorem \ref{theostable-constant}, it is sometimes
useful to represent the limiting L\'evy
motion $ ( \Lambda(\phi,u),   u\geq0 )$ in the form
%
\begin{equation} \label{eqlimit-processold}
\Lambda(\phi,u) = \int_0^u \int_\mcw\{\mce(w,\phi) -
\mce(0,\phi)\}   M_\alpha(\rmd s , \rmd w)  , \qquad   u \geq0   ,
\end{equation}
where $M_\alpha$ is a totally skewed to the right $\alpha$-stable
random measure on $(0,\infty)\times\mcw$ with control measure
$\lambda c_\alpha\mathrm{Leb}\times G$; see
Samorodnitsky and Taqqu \cite{samorodnitskytaqqu1994}. The representation~\eqref{eqlimit-processold} is linear in $\phi$, and this allows,
for example, handling more than one function $\phi$ at a time.\vadjust{\goodbreak}

Specifically, if Assumption \ref{hypoisp} holds, and $\mcf$ is a
class of
bounded measurable functions satisfying Assumption \ref{hypocty},
then, by
linearity, Theorem \ref{theostable-constant} implies that, for any
$n\geq2$
and bounded measurable functions $\phi_1,\ldots,\phi_n$ on $\mathcal
{D}([0,h])$
satisfying Assumption \ref{hypocty}, the family of $\bbr^n$-valued processes
$ (\mathcal Z_T(\phi_1,\cdot),\ldots, \mathcal Z_T(\phi
_n,\cdot) )$
converges weakly to the process $ (\Lambda(\phi_1,\cdot),\ldots,
\Lambda(\phi_n,\cdot) )$ in the sense of finite-dimensional
distributions. The components of the limiting process are defined by
\eqref{eqlimit-processold} and is an $\bbr^n$-valued $\alpha
$-stable L\'evy
motion. By Whitt \cite{whitt2002}, Theorem~11.6.7, the convergence also
holds in
$\mathcal D([0,\infty))^n$ endowed with the product (or weak) $M_1$ topology.

For another application of \eqref{eqlimit-processold}, we can write
the one-dimensional weak convergence prescribed by Theorem
\ref{theostable-constant} at $u=1$ in the form
%
\begin{equation} \label{eqlimit-process1d}
\mathcal Z_T(\phi,1) \Rightarrow\Lambda_1(\phi) :=\int_\mcw\{\mce
(w,\phi) -
\mce(0,\phi)\}   \tilde M_\alpha( \rmd w) ,\vspace*{-2pt}
\end{equation}
where this time $\tilde M_\alpha$ is a totally skewed to the right
$\alpha$-stable random measure on $\mcw$ with control measure
$\lambda c_\alpha  G$. Again, the representation of the limit in the
right-hand side of~\eqref{eqlimit-process1d} is linear in $\phi$, allowing us to handle
more than one function $\phi$ at a time.\vspace*{-4pt}
\end{remark}

\section{An application: The empirical process}\vspace*{-4pt}

Suppose we want to estimate the distribution function $K$ of
$X(0)$. For this purpose, we consider the family of empirical processes
\[
E_T(x) = T^{-1} \int_0^T \1_{\{ X(s) \leq x\}} \, \rmd s, \qquad    x>0 .\vspace*{-2pt}
\]
Let $D$
denote the set of discontinuity points of the distribution function
$K$ restricted to $\mcw\cap[0,\infty)$. The following is an immediate
corollary of Theorem \ref{theostable-constant} and
\eqref{eqlimit-process1d}.\vspace*{-2pt}
\begin{corollary} \label{corempirical}
Let $\mcx$ be the collection of $x>0$ such that $G(x-D)=0$. Then
\[
 \bigl( T a(T)^{-1}  \bigl( E_T(x)-K(x) \bigr),   x\in\mcx \bigr)
\Rightarrow \bigl( \mathcal{D}(x),   x\in\mcx \bigr)\vspace*{-2pt}
\]
in the sense of convergence of the finite-dimensional distributions,
where
\[
\mathcal{D}(x)=\int_\mcw\{K(x-w)
-K(x)\}\tilde M_\alpha(\rmd w), \qquad   x>0 .\vspace*{-2pt}
\]
\end{corollary}
\begin{remark}
Let us briefly comment on the condition $G(x-D)=0$.
\begin{enumerate}[3.]
\item[1.] Note that the set $D$ is at most countable, and the set of atoms of
$G$ is at most countable as well. We immediately conclude that the set
$\mcx$ misses at most countably many $x>0$.
\item[2.] Further, if the distribution of $W$ is supported on a closed set
consisting of isolated points in $[0,\infty)$, we have $D=\emptyset$ (see
Remark~\ref{remdiscCase}), and so $\mcx=(0,\infty)$.
\item[3.] Finally, $X(0)$ is an infinitely divisible random variable with
L\'evy measure
$\mu$ satisfying
\[
\mu ( (a,\infty) ) = \lambda\esp \bigl( Y\1(W>a) \bigr),
 \qquad
a>0 .\vspace*{-2pt}\vadjust{\goodbreak}
\]
Therefore, if $W$ does not have positive atoms, then the distribution function
$K$ has a single atom, at the origin, implying that $D=\{0\}$ and $\mcx
$ misses
some of the atoms of $G$, specifically those atoms that are not
isolated points
of $\mcw$.
\end{enumerate}
\end{remark}

\begin{remark}
It is important to note that estimators based on the empirical process $E_T$
may not be able to identify the parameter of interest, even for simple
parametric models of the distribution of $(Y,W)$. For instance,
if $Y$ and $W$ are independent, $K$ depends on the distribution of $Y$
only through its mean $\mathbb{E}[Y]$. This is the main motivation for
considering
the case $h>0$ in Theorem~\ref{theostable-constant} although it is
not the
object of this paper to provide practical details on this application.
\end{remark}

Observe that Corollary \ref{corempirical} shows that ``the usual''
$\sqrt{T}$-rate of convergence of an empirical process does not hold
in the
present situation, since the actual rate of convergence is~$Ta(T)^{-1}$, which
is regularly varying with index $1-\alpha^{-1}\in(0,1/2)$. This
should not be
surprising since presence of long range dependence has long been known
to yield
slower rates of convergence and non standard limit for the empirical
process. See, for example, Dehling and Taqqu \cite{dehlingtaqqu1989} for subordinated
Gaussian processes
and Surgailis \cite{surgailis2002,surgailis2004} for bounded functionals of
infinite or
finite variance linear processes.

\section{Proofs}\label{secproof-lemmas}

\begin{pf*}{Proof of Lemma~\ref{lemtailC}}
By the definition of $a$ and regular variation of the tail of $F$,
\[
\bar{F}(a(t)) = \pr\bigl(Y > a(t)\bigr) \sim t^{-1}    \qquad \mbox{as $t\to\infty$;}
\]
recall, further, that $a$ is regularly varying at infinity with index
$1/\alpha$. We will use the notation $N_{\epsilon,t}$,
$Y_{\epsilon,t}$, $\Gamma_{\epsilon,t}$ and $W_{\epsilon,t}$
introduced just before Lemma \ref{lemM1} above.
Applying Lemma~1 in Resnick and Samorodnitsky \cite{resnicksamorodnitsky1999a} and the
regular variation of $\bar F$, we get
%
\begin{equation} \label{eqNepsProba}
\lim_{t\to\infty} t\pr(N_{\epsilon,t} \geq1) = \lim_{t\to\infty}
\frac{\pr(N_{\epsilon,t} \geq1)}{\bar F(\epsilon a(t))} \frac{\bar
F(\epsilon a(t))}{\bar F(a(t))} = \rme^{\lambda\esp[Y]}
\epsilon^{-\alpha}   .
\end{equation}
Imagine, for a moment, that all sessions of the length exceeding
$\epsilon a(t)$ are
discarded upon arrival, and do not
contribute to a busy period.
Let $B_{\epsilon,t}$ denote the length of the first busy period
starting at or after
time $S_0$ and generated by the remaining
sessions, those of length not exceeding $\epsilon a(t)$. Then by
Resnick and Samorodnitsky \cite{resnicksamorodnitsky1999a}, Proposition~1, there exists a constant
$D$ independent of $\epsilon$ such that
%
\begin{equation}
\label{eqBepsBound}
\pr \bigl(B_{\epsilon,t}>\epsilon D a(t) \bigr)=\mathrm{o}(t^{-1})   .
\end{equation}
We immediately conclude that
%
\begin{equation}
\label{eqC1Nnul}
\lim_{t\to\infty} t\pr\bigl(C_1>\epsilon D a(t)   ; N_{\epsilon,t} =
0\bigr) = 0
\end{equation}
(keeping in mind that an idle period has an exponential
distribution).\vadjust{\goodbreak}

We consider now the case $N_{\epsilon,t}\geq1$, in which case we use
the decomposition
%
\begin{equation}
\label{eqC1decomp}
C_1=\{\Gamma_{\epsilon,t}-S_0\} + Y_{\epsilon,t} +
\{S_1-(\Gamma_{\epsilon,t}+Y_{\epsilon,t}) \}  .
\end{equation}

Since $B_{\epsilon,t}$ is the length of the first busy session
starting after
$S_0$ and generated only by sessions of length less than $\epsilon
a(t)$ and
since $\Gamma_{\epsilon,t}$ is the starting point of the first
session of
length greater than $\epsilon a(t)$ starting after $S_0$, it is clear that
$S_0+B_{\epsilon,t} < \Gamma_{\epsilon,t} $ implies that
$N_{\epsilon,t}=0$.
Thus, on the event $\{N_{\epsilon,t}\geq1\}$, it holds that
\[
\Gamma_{\epsilon,t}-S_0\leq B_{\epsilon,t} .
\]
Hence, by (\ref{eqBepsBound}), for any $\eta>0$, choosing
$\epsilon>0$ sufficiently small (i.e., $\epsilon<\eta/D$ where $D$ is
as in (\ref{eqC1Nnul})), we have
%
\begin{equation} \label{estep22}
\pr \bigl(\Gamma_{\epsilon,t}-S_0 >a(t)\eta;N_{\epsilon,t}\geq
1 \bigr) =
\mathrm{o}(t^{-1})\qquad\mbox{as $t\to\infty$} .
\end{equation}
Further, denote by $\tilde{\Gamma}_{\epsilon,t}$ the completion time
of the last session with length greater than~$\epsilon a(t)$ before
time $S_1$. Notice that the infinite source Poisson process
\eqref{eqdefgginfini} is time reversible, in the sense of switching
the direction of time, declaring $\Gamma_\ell+ Y_\ell$ to be the
arrival time of session number $\ell$ and $\Gamma_\ell$ to be its
completion time. Therefore, by time inversion, the difference
$S_1-\tilde{\Gamma}_{\epsilon,t}$ has the same distribution
as
$\Gamma_{\epsilon,t}-S_0+I_0$, where $I_0$ denotes the idle period
preceding $S_0$. Moreover, the joint distribution of
$(S_1-\tilde{\Gamma}_{\epsilon,t}, N_{\epsilon,t})$
and $(\Gamma_{\epsilon,t}-S_0+I_0, N_{\epsilon,t})$ are also the same.
Since on the event $\{ N_{\epsilon,t} = 1\}$, the
random variables $\Gamma_{\epsilon,t}+Y_{\epsilon,t}$ and
$\tilde{\Gamma}_{\epsilon,t}$ coincide, we conclude that, for all
$\eta,\epsilon>0$,
%
\begin{eqnarray}
\label{estep23}
&&\pr \bigl(S_1 - (\Gamma_{\epsilon,t} + Y_{\epsilon,t}) > a(t) \eta
  ;
N_{\epsilon,t} = 1  \bigr)\nonumber\\
&& \quad = \pr \bigl( S_1 - \tilde{\Gamma}_{\epsilon,t} > a(t) \eta;
N_{\epsilon,t} =
1  \bigr) \nonumber
\\[-8pt]
\\[-8pt]
&& \quad = \pr \bigl( \Gamma_{\epsilon,t} - S_0+I_0 > a(t) \eta; N_{\epsilon
,t} =
1  \bigr) \nonumber\\
&& \quad \leq
\pr \bigl( \Gamma_{\epsilon,t} - S_0 > a(t) \eta/2; N_{\epsilon,t}
\geq1
 \bigr) + \pr \bigl( I_0 > a(t) \eta/2
 \bigr)
= \mathrm{o}(t^{-1}) \qquad\mbox{as $t\to\infty$} ,
\nonumber
\end{eqnarray}
where the $\mathrm{o}$-term follows from (\ref{estep22}) and the fact that
$I_0$ has
exponential distribution. Next, by Lemma~2 in
Resnick and Samorodnitsky \cite{resnicksamorodnitsky1999a}, we also have
%
\begin{equation}\label{eqmoinsde2}
\pr(N_{\epsilon,t} \geq2) = \mathrm{o}(t^{-1})\qquad\mbox{as $t\to\infty
$} .
\end{equation}
Applying (\ref{eqC1Nnul}), (\ref{eqC1decomp}), (\ref
{estep22}), (\ref{estep23})
and (\ref{eqmoinsde2}),
we get, for any $x>\eta>0$, choosing $\epsilon$ small enough,
%
\begin{eqnarray}
\label{eqintermediaire-B}
\liminf_{t\to\infty} t \pr\bigl(Y_{\epsilon,t} > a(t) x   ;
N_{\epsilon,t} \geq1\bigr)
&\leq&\liminf_{t\to\infty} t  \pr\bigl(C_1>a(t)x\bigr) \nonumber
\\
&\leq&\limsup_{t\to\infty} t  \pr\bigl(C_1>a(t)x\bigr)\\
 &\leq&
\limsup_{t\to\infty} t  \pr\bigl(Y_{\epsilon,t} > a(t) (x-\eta)   ;
N_{\epsilon,t} \geq1\bigr)   .
\nonumber
\end{eqnarray}
Note that the distribution of $Y_{\epsilon,t}$ is the conditional distribution
of $Y$ given $\{Y>\epsilon a(t)\}$ and that the event $\{N_{\epsilon,t}
\geq1\}$ is independent of $Y_{\epsilon,t}$, so that (\ref{eqNepsProba})
yields, for any $x>0$,
\begin{eqnarray*}
t  \pr\bigl(Y_{\epsilon,t} > a(t) x   ;   N_{\epsilon,t} \geq1\bigr) \sim
\rme^{\lambda\esp[Y]}\epsilon^{-\alpha} \pr\bigl(Y >a(t) x |
Y>\epsilon a(t) \bigr) \to\rme^{\lambda\esp[Y]} x^{-\alpha}\vspace*{-2pt}
\end{eqnarray*}
as $t\to\infty$.
Applying this statement to
 (\ref{eqintermediaire-B}) and letting $\eta\to0$ gives (\ref{ecyclel}).\vspace*{-2pt}
\end{pf*}

\begin{pf*}{Proof of Lemma~\ref{lemtruly-iid}}
Observe that the process $\{X(t),t\in\Rset\}$ is a regenerative
process (it
regenerates at the beginning of each busy period), hence it is
ergodic. Therefore, $T^{-1} \mathcal{J}_T(\phi) \to\mathcal E(0,\phi
)$ a.s.;
see, for example, Resnick \cite{resnick1992}. On the other hand, as seen
earlier, the
sequence $ (Z_j(\phi) )$ is strongly mixing, hence also
ergodic, and
so $n^{-1}\sum_{j=1}^{n} Z_j(\phi)$ converges almost surely to
$\esp[Z_1(\phi)]$. For $T>0$, let $M_T$ denote the number of complete cycles
initiated after time $0$, and finishing before time $T$. Since $M_T/T$
converges almost surely to $1/\esp[C_1]$, we also obtain
\[
\frac1T \sum_{j=1}^{M_T} Z_j(\phi) \to\esp[Z_1(\phi)]/\esp[C_1]
 , \qquad
 \mbox{a.s.,}\vspace*{-2pt}
\]
and (\ref{eqZmean}) follows.

Denote $\bar\phi= \phi- \mathcal E(0,\phi)$. Observe that $\mathcal
{J}_T(\bar
\phi)$ is centered and
$\|\bar\phi\|_\infty\leq\|\phi\|_\infty+|\mce(0,\phi)|\leq2\|
\phi\|_\infty$.
We have
%
\begin{equation}
\label{eqRT1}
\sup_{t\in[0,S_0]} |\mathcal{J}_{t}(\bar\phi) |\leq
S_0\|\bar\phi\|_\infty .\vspace*{-2pt}
\end{equation}
For $t\geq S_0$, we use the decomposition
\[
\mathcal{J}_{t}(\bar\phi) = \mathcal{J}_{S_0}(\bar\phi) + \sum
_{j=1}^{M_{t}}
Z_j(\bar\phi)+\int_{S_{M_{t}}}^{t}\bar\phi(X_h(s))\,\rmd s .\vspace*{-2pt}
\]
Now, using $\|\bar\phi\|_\infty\leq2\| \phi\|_\infty$, (\ref{eqRT1})
and that, for all $k=1,\ldots,M_T+1$,
\[
\sup_{u\in[S_{k-1},S_k]} \biggl|\int_{S_{k-1}}^{u}\bar\phi
(X_h(s))\,\rmd s \biggr|\leq
\|\bar\phi\|_\infty C_k  ,\vspace*{-2pt}
\]
we get, for any $T>0$,
\begin{eqnarray*}
\pr \Bigl(\sup_{t\in[0,T]} |\mathcal{J}_{t}(\bar\phi)
| > 5 x
\|\phi\|_\infty \Bigr)
& \leq&\pr(S_0 > x) + \pr \Biggl(\sup_{t\in[0,T]}  \Biggl|\sum
_{j=1}^{M_{t}} Z_j(\bar
\phi) \Biggr| >
x \|\phi\|_\infty \Biggr) \\[-2pt]
&&{}  + \pr\Bigl (\max_{k=1,\ldots,M_T+1}C_k> x \Bigr)\\[-2pt]
& \leq& \pr(S_0 > x) + 2 \pr(M_T > 2T/\esp[C_1]) \\[-2pt]
&&{}   + \pr \Biggl( \max_{1\leq k \leq2T/\esp[C_1]
} \Biggl| \sum_{j=1}^{k}
Z_j(\bar\phi)  \Biggr| > x \|\phi\|_\infty \Biggr) \\[-2pt]
&&{}  + (2T/\esp[C_1]+1)\pr (C_1> x )  .\vspace*{-2pt}\vadjust{\goodbreak}
\end{eqnarray*}
Applying (\ref{eqZmean}), we see that $Z_j(\bar\phi)$ is centered.
Moreover, $|Z_j(\bar\phi)|\leq2C_j\|\phi\|_\infty$. Let
$p\in(1,\alpha)$. Applying the mixing property (\ref{eqmixing}),
Lemma~\ref{lemtailC} and Rio \cite{rio2000}, Chapiter~3, Exercise~1,
there exists a constant $c$ which depends only on the distribution of
$C_1$ and $p$ such that
%
\begin{equation} \label{eqburkholder-mixing}
\esp \Biggl[ \max_{1\leq k \leq n }  \Biggl| \sum_{j=1}^{k} Z_j(\bar
\phi)  \Biggr|^p  \Biggr] \leq c \|\phi\|_\infty^p n   .
\end{equation}
Finally, we bound $\pr(M_T > 2T/\esp[C_1])$ by noting as usual that
$M_T>n$ if and only if $S_{n+1} \leq T$. Thus, denoting by $m$ the
smallest integer larger than or equal to $2T/\esp[C_1]$, we have, for
some constant $c$ only depending on the distribution of $C_1$ and $p$,
\[
\pr(M_T>2T/\esp[C_1]) \leq\pr(S_m \leq T) \leq\pr(S_m - m\esp[C_1]
\leq-T) \leq T^{-p}\esp\bigl[|S_m-m\esp[C_1]|^p\bigr] .
\]
Since $S_m-m\esp[C_1]$ is a sum of i.i.d. centered random variables
with finite
$p$th moment, we obtain by Burkh\"older inequality
(see~von Bahr and Esseen \cite{vonbahresseen1965}, Theorem~2),
%
\begin{equation}\label{eqMTbound}
\pr(M_T>2T/\esp[C_1]) =\mathrm{O}(T^{1-p}) .
\end{equation}
Gathering the previous displays and using
$\pr (C_1> x )\leq\esp[C_1^p]x^{-p}$ for any $p<\alpha$, we
obtain~(\ref{equnifBound}) with $g(x) = P(S_0>x)$.
\end{pf*}
\begin{pf*}{Proof of Lemma~\ref{lemM1}}
We will bound the function
\[
\Delta(v)=\int_{S_{0}}^{v}\bigl \{\phi(X_h(s)) -
\mce(W_{\epsilon,t},\phi) \1_{[\Gamma_{\epsilon,t},
\Gamma_{\epsilon,t} + Y_{\epsilon,t})}(s)\bigr\}\, \rmd s
\]
on the event $\{N_{\epsilon,t}\geq1\}$ successively for
$v\in[S_0,\Gamma_{\epsilon,t}]$,
$v\in[\Gamma_{\epsilon,t},\Gamma_{\epsilon,t}+Y_{\epsilon,t}]$ and
$v\in[\Gamma_{\epsilon,t}+Y_{\epsilon,t},S_1]$.

\textit{Step} 1. For $v\in[S_0,\Gamma_{\epsilon,t}]$, we have
\[
|\Delta(v)|= \biggl|\int_{S_{0}}^{v} \phi(X_h(s)) \,\rmd s \biggr|\leq
(\Gamma_{\epsilon,t}-S_0)\|\phi\|_\infty .
\]
Hence, using (\ref{estep22}), for any $\eta>0$, choosing $\epsilon
>0$ sufficiently
small, we have
%
\begin{equation} \label{eqsupFirstIntervalResult}
\pr \Bigl(\sup_{v\in[S_0,\Gamma_{\epsilon,t}]}|\Delta(v)|>
a(t)\eta;N_{\epsilon,t}\geq1 \Bigr) = \mathrm{o}(t^{-1}) .
\end{equation}

\textit{Step} 2. For $v\in[\Gamma_{\epsilon,t}, \Gamma_{\epsilon,t}+
Y_{\epsilon,t}]$,
we write
%
\begin{eqnarray}
\label{eqstep20}
|\Delta(v)|&\leq&|\Delta(\Gamma_{\epsilon,t})|+|\Delta(v)-\Delta
(\Gamma_{\epsilon,t})|\nonumber
\\[-8pt]
\\[-8pt]
&\leq&\sup_{v\in[S_0,\Gamma_{\epsilon,t}]}|\Delta(v)| +
\sup_{y\in[0, Y_{\epsilon,t}]} \biggl|\int_{0}^{y}\bigl\{\phi\bigl(X_h(\Gamma
_{\epsilon,t}+s)\bigr) -
\mce(W_{\epsilon,t},\phi)\bigr\} \, \rmd s \biggr| .
\nonumber
\end{eqnarray}
For $s\in(0,Y_{\epsilon,t})$, $X(\Gamma_{\epsilon,t}+s)$ can be
expressed as
\[
X(\Gamma_{\epsilon,t}+s) = W_{\epsilon,t} + \check X(s) + R(s)   ,
\]
where $R(s)$
is the sum of all transmission rates of the sessions that started
before time~$\Gamma_{\epsilon,t}$ and are still active at time $s$,
and $\{\check X(s), s\geq0\}$ is defined by
\[
\check X(s) =
\sum_{\ell\in\Zset}W_\ell\1_{\{\Gamma_{\epsilon,t}<\Gamma_\ell
\leq
s+\Gamma_{\epsilon,t}<\Gamma_\ell+Y_\ell\}}   .
\]
Since each session that arrives after time $S_0$ but before time
$\Gamma_{\epsilon,t}$ has a length not exceeding $\epsilon a(t)$, we
conclude that $R(s)=0$ for $s> \epsilon a(t)$.
Using the notation $\check X_h(s) = \{\check X(s+v)   , 0 \leq v \leq
h\}$, we, therefore, obtain
%
\begin{equation}\label{eqstep21}
\sup_{y\in[0, Y_{\epsilon,t}]}
 \biggl| \int_{0}^{y}
\bigl\{\phi\bigl(X_h(\Gamma_{\epsilon,t}+s)\bigr) - \phi\bigl(W_{\epsilon,t} + \check
X_h(s) \bigr)\bigr\}
\, \rmd s  \biggr| \leq2 \|\phi\|_\infty
\epsilon a(t)   .
\end{equation}
Observe that the process $\check X$ is independent of
$(Y_{\epsilon,t}, W_{\epsilon,t}, \1_{\{N_{\epsilon,t}\geq1\}})$. We
preserve this independence
while transforming $\check X$ into a stationary process, with the same
law as the original process $X$ in \eqref{eqdefgginfini} by defining
\[
\hat X(s) = \sum_{\ell\leq0} W'_\ell  \1_{\{
\Gamma'_\ell\leq s < \Gamma'_\ell+ Y'_\ell\}} + \check X(s)
,\qquad s\in\Rset ,
\]
where $\{(\Gamma'_\ell,Y'_\ell,W'_\ell), \ell\in\Zset\}$ is an
independent copy of $\{(\Gamma_\ell,Y_\ell,W_\ell), \ell\in\Zset
\}$.
Clearly,
\[
\sup_{y\in[0,Y_{\epsilon,t}]} \biggl|\int_{0}^{y} \bigl\{\phi
\bigl(W_{\epsilon,t} + \check
X_h(s)\bigr)-\phi\bigl(W_{\epsilon,t} + \hat X_h(s)\bigr)\bigr\}\, \rmd s \biggr| \leq2
\|\phi\|_\infty\sup_{\ell\leq0}(\Gamma'_\ell+Y'_\ell)_+ ,
\]
where $\hat X_h(s) = \{\hat X(s+v)   , 0 \leq v \leq h\}$.
The random variable in the right-hand side above is finite with
probability 1 and independent of $N_{\epsilon,t}$. Therefore, it
follows from (\ref{eqNepsProba}) that for any $u>0$,
\[
\pr \Bigl( \sup_{\ell\leq0}(\Gamma'_\ell+Y'_\ell)> a(t) u   ;
N_{\epsilon,t} \geq1 \Bigr)= \mathrm{o}(t^{-1}) .
\]
The last two displays and (\ref{eqstep21}) give that, for any $\eta
>0$ and
$0<\epsilon<\eta/(2\|\phi\|_\infty)$,
%
\begin{eqnarray}\label{eqstep22}
&&\pr \biggl(\sup_{y\in[0, Y_{\epsilon,t}]}
 \biggl| \int_{0}^{y}
\bigl\{\phi\bigl(X_h(\Gamma_{\epsilon,t}+s)\bigr) - \phi\bigl(W_{\epsilon,t} + \hat
X_h(s) \bigr)\bigr\}
\,\rmd s  \biggr| >a(t)\eta
;N_{\epsilon,t}\geq1 \biggr)\nonumber
\\[-8pt]
\\[-8pt]
&& \quad =\mathrm{o}(t^{-1})   .\qquad
\nonumber
\end{eqnarray}

The event $\{N_{\epsilon,t}\geq1\}$ is, clearly,
independent of $(Y_{\epsilon,t},W_{\epsilon,t})$. Furthermore, the
latter pair has the conditional distribution of $(Y,W)$
given that $\{Y>\epsilon a(t)\}$. Since
$\hat X$ has the same law as $X$, we get for any $x>0$,
%
\begin{eqnarray}\label{eqStep23}
&&\pr \biggl(\sup_{y\in[0, Y_{\epsilon,t}]}  \biggl|\int_{0}^{y}
 \bigl\{\phi\bigl(W_{\epsilon,t} + \hat X_h(s)\bigr) -
\mce(W_{\epsilon,t},\phi)  \bigr\} \, \rmd s  \biggr| > x  ;
N_{\epsilon,t} \geq1 \biggr)
\nonumber
\\[-8pt]
\\[-8pt]
&& \quad = \pr \biggl(\sup_{y\in[0, Y]} \biggl| \int_0^y  \bigl\{ \phi\bigl(W +
X_h(s)\bigr)- \mce(W,\phi)  \bigr\} \, \rmd s  \biggr| > x \bigm| Y>
\epsilon a(t)  \biggr) \times\pr(N_{\epsilon,t} \geq1)  ,
\nonumber
\end{eqnarray}
where the pair $(Y,W)$ in the right-hand side is taken to be
independent of the\break process~$X$.

Recall that $\mce(w,\phi) = \esp[\phi(w+X_h(0))]$, that for any
$w\geq0$, $\|\phi(w+\cdot)\|_\infty\leq\|\phi\|_\infty$ and, for any
$y\geq0$, $\esp[\mathcal{J}_y(\phi(w+\cdot))]=y\mce(w,\phi)$. It
follows from these observations and (\ref{equnifBound}) in
Lemma~\ref{lemtruly-iid} that, for any $x>0$,
\[
\sup_{w\geq0}\pr \biggl(\sup_{y\in[0,u]} \biggl|\int_0^y \bigl\{\phi\bigl(w+
X_h(s)\bigr)-\mce(w,\phi)\bigr\}\, \rmd s
 \biggr|>x\|\phi\|_\infty \biggr)\leq C u^{1-p}+C u x^{-p}+g(x) ,
\]
for $p\in(1,\alpha)$, some constant $C>0$ and $g(x)\to0$ as
$x\to\infty$. Integrating in $(w,u)$ with respect to the distribution
of $(W,Y)$ in \eqref{eqStep23}, this bound yields, for any $u>0$ and
$A>0$,
\begin{eqnarray*}
&&\pr \biggl( \sup_{y\in[0, Y]}  \biggl| \int_0^y \bigl\{\phi\bigl(W + X_h(s)\bigr) -
\mce(W,\phi)\bigr\}\, \rmd s  \biggr|> u A \bigm| Y> A  \biggr)
\\
&& \quad \leq C \esp[Y^{1-p} \mid Y>A]+C \|\phi\|_\infty^p
(uA)^{-p}\esp[Y \mid Y>A] + g(uA/\|\phi\|_\infty)  .
\end{eqnarray*}
As $A\to\infty$, we have both $\esp[Y^{1-p} \mid Y>A]\to0$ and
$A^{-p}\esp[Y \mid Y>A]\to0$ since $Y$ has a~regularly varying tail
with index $\alpha>1$ and $p\in(1,\alpha)$. Thus, the 3 terms in the
previous bound converge to 0 as $A\to\infty$. This, together
with (\ref{eqStep23}) and (\ref{eqNepsProba}), yields that, for any
$\epsilon>0$ and $\eta>0$,
\[
\pr \biggl(\sup_{y\in[0, Y_{\epsilon,t}]}
 \biggl|\int_{0}^{y}  \bigl\{\phi\bigl(W_{\epsilon,t} + \hat X_h(s)\bigr) -
\mce(W_{\epsilon,t},\phi)  \bigr\} \, \rmd s  \biggr| > a(t) \eta
  ;
N_{\epsilon,t} \geq1 \biggr)=\mathrm{o}(t^{-1}) .
\]
Finally, gathering the last display, (\ref{eqstep22}), (\ref{eqstep20})
and (\ref{eqsupFirstIntervalResult}), we obtain
%
\begin{equation} \label{eqsupSecIntervalResult}
\pr \Bigl(\sup_{v\in[\Gamma_{\epsilon,t},\Gamma_{\epsilon
,t}+Y_{\epsilon,t}]}|\Delta(v)|>
a(t)\eta;N_{\epsilon,t}\geq1 \Bigr) = \mathrm{o}(t^{-1}) .
\end{equation}

\textit{Step} 3. If $v\in[\Gamma_{\epsilon,t}+Y_{\epsilon,t},S_1]$, we
have on
$\{N_{\epsilon,t}\geq1\}$,
%
\begin{eqnarray}\label{eqlastIntervalSup}
|\Delta(v)|&\leq&|\Delta(\Gamma_{\epsilon,t}+Y_{\epsilon,t})| +
 \biggl|\int_{\Gamma_{\epsilon,t}+Y_{\epsilon,t}}^{v} \phi(X_h(s))
\,\rmd s \biggr|\nonumber
\\[-8pt]
\\[-8pt]
&\leq&\sup_{v\in[\Gamma_{\epsilon,t},\Gamma_{\epsilon
,t}+Y_{\epsilon,t}]}|\Delta(v)|
+\{S_1-(\Gamma_{\epsilon,t}+Y_{\epsilon,t})\} \|\phi\|_\infty .
\nonumber
\end{eqnarray}
Using (\ref{eqsupSecIntervalResult}) (\ref
{eqlastIntervalSup}), (\ref{estep23})
and (\ref{eqmoinsde2}), for any $\eta>0$, we have
%
\begin{equation} \label{eqsupLastIntervalResult}
\pr \Bigl(\sup_{v\in[\Gamma_{\epsilon,t}+Y_{\epsilon
,t},S_1]}|\Delta(v)|>
a(t)\eta;N_{\epsilon,t}\geq1 \Bigr) = \mathrm{o}(t^{-1}) .
\end{equation}
\end{pf*}
\begin{pf*}{Proof of Lemma~\ref{lemjoint-tails}}
Let $f$ a Lipschitz function with compact support in
$[-\infty,\infty]^d\setminus\{0\}$, and let $L$ be its Lipschitz
constant. Let $c>0$ be small enough such that the support of~$f$
does not intersect $[-2c,2c]^d$.

Using the fact that, in the notation of \eqref{eqZj},
$|Z_1(\phi_i)|\leq\|\phi_i\|_\infty
C_1$ for each $i=1,\ldots,d$, the bound (\ref{eqC1Nnul}) implies
that, as $t\to\infty$,
\[
\pr \bigl( |Z_1(\phi_i) |>c a(t)  \mbox{ for some
$i=1,\ldots, d$} ; N_{\epsilon,t}=0 \bigr)=\mathrm{o}(t^{-1})
\]
as long as $\epsilon>0$ is small enough relatively to $c$. We will
show that
%
\begin{eqnarray}
\label{eprovelemma4}
&&\lim_{\epsilon\to0}   \limsup_{t\to\infty} t \esp\bigl[ f\bigl(\mathbf Z /
a(t)\bigr)   ;
  N_{\epsilon,t} \geq1 \bigr] \nonumber\\
&& \quad  = \lim_{\epsilon\to0} \liminf_{t\to\infty} t \esp\bigl[f\bigl(\mathbf Z /
a(t)\bigr)   ;   N_{\epsilon,t} \geq1 \bigr] \\
&& \quad  = \rme^{\lambda\esp[Y]}   \int_0^\infty\esp [f(y[\mce
(W^*,\phi_1),\ldots, \mce(W^*,\phi_1)]^{\mathrm{T}}) ]\alpha y^{-\alpha-1}\,\rmd y
  . \nonumber
\end{eqnarray}
This will prove the required vague convergence in \eqref{eqqueueZmultdim}.
Write
%
\begin{eqnarray}\label{eqslutskyArgument}
t \esp \bigl[ f\bigl(\mathbf Z / a(t)\bigr) ; N_{\epsilon,t} \geq1  \bigr]
&=& t \esp \bigl[ f\bigl(\boldsymbol\Phi(Y_{\epsilon,t}, W_{\epsilon,t}) /
a(t)\bigr)   ;   N_{\epsilon,t} \geq1  \bigr] \nonumber
 \quad \\[-8pt]
\\[-8pt]
&&{}+ t \esp \bigl[\bigl \{ f\bigl(\mathbf Z / a(t)\bigr) - f\bigl(\boldsymbol\Phi(Y_{\epsilon,t}, W_{\epsilon,t}) / a(t)\bigr) \bigr\}   ;
N_{\epsilon,t} \geq1  \bigr]  , \quad
\nonumber
\end{eqnarray}
where $\boldsymbol{\Phi}(y,w) = y [
\mce(w,\phi_1),\ldots,\mce(w,\phi_d) ]^{\mathrm{T}}$. Choose $0<\eta<c$
and observe
that the Lipschitz property of $f$ and the fact that its support does not
intersect $[-2c,2c]^d$ implies that, on the event
$\bigcap_i\{|Z_1(\phi_i)-\mce(W_{\epsilon,t},\phi_i)Y_{\epsilon
,t}|\leq\eta
a(t)\}$,
\[
\bigl|f\bigl(\mathbf Z / a(t)\bigr) - f\bigl(\boldsymbol\Phi(Y_{\epsilon,t},W_{\epsilon,t}) / a(t)\bigr)\bigr| \leq L\eta  \1
 \bigl(|\mce(W_{\epsilon,t},\phi_i)Y_{\epsilon,t}| > \eta a(t)
\mbox{ for some $i=1,\ldots, d$} \bigr)   .
\]
Letting $g$ be a continuous function on $[-\infty,\infty]^d$
such that $g(x)=1$ for all $x\notin[-c,c]^d$
and $g(x)=0$ in a neighborhood of the origin, we obtain
\begin{eqnarray*}
&&t\esp \bigl[ \bigl|f\bigl(\mathbf Z / a(t)\bigr)-f\bigl(\boldsymbol\Phi(Y_{\epsilon,t},W_{\epsilon,t}) / a(t)\bigr)  \bigr|   ;
N_{\epsilon,t} \geq1  \bigr]\\
 && \quad \leq L \eta  t\esp \bigl[
g\bigl(\boldsymbol\Phi(Y_{\epsilon,t},W_{\epsilon,t}) / a(t)\bigr) ;
N_{\epsilon,t} \geq1  \bigr]  \\
&& \qquad {}+2 \|f\|_\infty\sum_{i=1}^d t  \pr \bigl( |Z_1(\phi_i) -
\mce(W_{\epsilon,t},\phi_i) Y_{\epsilon,t}  | > \eta a(t)
 ;   N_{\epsilon,t} \geq1  \bigr)   .
\end{eqnarray*}
Recall that by Lemma~\ref{lemM1},
\[
\lim_{t\to\infty} t \pr \bigl( |Z_1(\phi_i) -
\mce(W_{\epsilon,t},\phi_i) Y_{\epsilon,t}  | > \eta a(t)
  ;   N_{\epsilon,t} \geq1  \bigr) = 0
\]
for all $\epsilon>0$ small enough (relative to $\eta$). Therefore, for
each $\eta>0$ and $\epsilon>0$ small enough,
\begin{eqnarray*}
&&\limsup_{t\to\infty}  \bigl| t \esp \bigl[ f\bigl(\mathbf Z / a(t)\bigr)  ;
N_{\epsilon,t} \geq1  \bigr] - t \esp \bigl[ f\bigl(\boldsymbol\Phi(Y_{\epsilon,t}, W_{\epsilon,t}) / a(t)\bigr)   ;
N_{\epsilon,t} \geq1  \bigr]  \bigr|
\\
&& \quad \leq L\eta\limsup_{t\to\infty} t\esp \bigl[g\bigl(\boldsymbol\Phi(Y_{\epsilon,t},W_{\epsilon,t}) /
a(t)\bigr) ; N_{\epsilon,t}\geq1 \bigr] .
\end{eqnarray*}

We will prove below that for any $\epsilon>0$,
%
\begin{equation}
\label{eqcontmapping}
t \pr \biggl( \frac{\boldsymbol\Phi
(Y_{\epsilon,t},W_{\epsilon,t})}{a(t)} \in\cdot  ;
N_{\epsilon,t} \geq1  \biggr) \stackrel v \longrightarrow
\rme^{\lambda\esp[Y]}    (\nu_{\alpha;\epsilon}\times
G )
\circ\boldsymbol\Phi^{-1}(\cdot)  ,
\end{equation}
where the measure $\nu_{\alpha;\epsilon}$ on $(0,\infty)$ is the
restriction of the measure $\nu_\alpha$ in \eqref{eqdef-G} to
$(\epsilon,\infty)$, i.e.
$\nu_{\alpha;\epsilon}(x,\infty) = \min (x^{-\alpha},
\epsilon^{-\alpha} )$, $x>0$. Assuming this has been proved, it
will follow that
%
\begin{eqnarray}
\label{elem4g}
&&\limsup_{t\to\infty}  \bigl| t  \esp \bigl[ \bigl\{ f\bigl(\mathbf Z / a(t)\bigr) -
f\bigl(\boldsymbol\Phi(Y_{\epsilon,t}, W_{\epsilon,t}) / a(t)\bigr) \bigr\}   ;
N_{\epsilon,t} \geq1  \bigr]  \bigr|\nonumber\\
&& \quad  \leq CL\eta \int g\circ\boldsymbol\Phi\, \rmd(\nu_{\alpha
;\epsilon}\times G) \\
&& \quad \leq CL\eta \int g\circ\boldsymbol\Phi\, \rmd(\nu_\alpha\times G)
\nonumber
\end{eqnarray}
for some finite positive constant $C$ independent of $\eta$ and
$\epsilon$. Note that the last integral is finite. Similarly,
\eqref{eqcontmapping} will imply that
%
\begin{eqnarray} \label{elem4f}
\lim_{t\to\infty} t  \esp \bigl[ f\bigl(\boldsymbol\Phi
(Y_{\epsilon,t},W_{\epsilon,t}) / a(t)\bigr)   ;
N_{\epsilon,t} \geq1  \bigr] &=& \rme^{\lambda\esp[Y]}   \int
f \circ\boldsymbol{\Phi} \,\rmd(\nu_{\alpha;\epsilon}\times G)
\nonumber
\\[-8pt]
\\[-8pt]
&=& \rme^{\lambda\esp[Y]}   \int
f \circ\boldsymbol{\Phi} \,\rmd(\nu_{\alpha}\times G)
\nonumber
\end{eqnarray}
for all $0<\epsilon<c/(\max_{i=1,\ldots,d}\|\phi_i\|_\infty)$. We combine
\eqref{eqslutskyArgument}, \eqref{elem4g} and \eqref{elem4f} by
keeping $\eta$ fixed and letting $\epsilon\to0$. This shows that
\begin{eqnarray*}
&&-C  L\eta\int g\circ\boldsymbol\Phi\, \rmd(\nu_\alpha\times G)
+ \rme^{\lambda\esp[Y]}   \int
f \circ\boldsymbol{\Phi} \,\rmd(\nu_{\alpha}\times G)\\
&& \quad \leq\lim_{\epsilon\to0}\liminf_{t\to\infty}
t\esp \bigl[f\bigl(\mathbf Z / a(t)\bigr) ; N_{\epsilon,t}\geq1 \bigr]
\\
&& \quad \leq\lim_{\epsilon\to0}\limsup_{t\to\infty} t\esp
\bigl[f\bigl(\mathbf Z /
a(t)\bigr) ; N_{\epsilon,t}\geq1 \bigr]\\
&& \quad \leq C  L\eta\int g\circ\boldsymbol\Phi\, \rmd(\nu_\alpha\times G)
+ \rme^{\lambda\esp[Y]}   \int
f \circ\boldsymbol{\Phi} \,\rmd(\nu_{\alpha}\times G) ,
\end{eqnarray*}
and \eqref{eprovelemma4} follows by letting $\eta\to0$.

It remains to prove  (\ref{eqcontmapping}).
holds. Since the event $\{N_{\epsilon,t}\geq1\}$ is
independent of $(Y_{\epsilon,t},W_{\epsilon,t})$, whose distribution
is the conditional distribution of $(Y,W)$ given that $\{Y>\epsilon
a(t)\}$, we have, as $t\to\infty$,
\begin{eqnarray*}
t\pr \bigl(\boldsymbol\Phi\bigl(Y_{\epsilon,t}/a(t),W_{\epsilon,t}\bigr)\in\cdot ; N_{\epsilon
,t}\geq1 \bigr)
&=&
t \pr (N_{\epsilon,t}\geq1 )\times\pr \bigl(\boldsymbol\Phi\bigl(Y/a(t),W\bigr)\in\cdot
\mid Y>\epsilon a(t) \bigr)\\
&\sim&\rme^{\lambda\esp[Y]}
\epsilon^{-\alpha} \pr\bigl (\boldsymbol\Phi\bigl(Y/a(t),W\bigr)\in\cdot\mid Y>\epsilon a(t) \bigr) ,
\end{eqnarray*}
by (\ref{eqNepsProba}). Further, by
Assumption~\ref{hypoisp}\ref{itemjointe},
\[
\pr \bigl(\bigl(Y/a(t),W\bigr)\in\cdot\mid Y>\epsilon a(t) \bigr) \stackrel v
\longrightarrow\epsilon^{\alpha} \nu_{\alpha;\epsilon}\times G .
\]

We extend $\boldsymbol\Phi$ to
$(0,\infty)\times[0,\infty]$ by
\[
\boldsymbol\Phi(y,\infty)=\lim_{w\to\infty}\boldsymbol\Phi(y,w) ,
\]
when the limit exists, or by defining the value at infinity to be equal
to 0
otherwise. Then the set of discontinuities of $\boldsymbol\Phi$ in
$(0,\infty]\times\mcw$ is included in
\[
(0,\infty)\times\bigcup_{i=1,\ldots,d}D(\mce(\cdot,\phi_i),\mcw
) ,
\]
which has $\nu_{\alpha;\epsilon}\times G$-measure zero by (\ref{eqdef-G}),
since each function $\phi_i$ satisfies Assumption~\ref{hypocty}.
Now, since
$\boldsymbol\Phi(y,w)/a(t)=\boldsymbol\Phi(y/a(t),w)$, (\ref
{eqcontmapping}) follows
from the continuous mapping theorem.
\end{pf*}

\begin{pf*}{Proof of Theorem~\ref{theostable-constant}}
In order to prove convergence in $\mathcal D([0,\infty))$ it is enough to
prove convergence in $\mathcal D([0,a])$ for any $a>0$. For notational
simplicity, we present the argument for $a=1$.

For any bounded interval $[a,b]$ and real-valued functions $x_1$ and
$x_2$ in
$\mathcal D([a,b])$, we denote by $d_{M_1}(x_1,x_2,[a,b])$ the $M_1$ distance
between $x_1$ and $x_2$ on $[a,b]$, and we write $d_{M_1}(x_1,x_2)$ if
$[a,b]=[0,1]$. We refer the reader to Whitt \cite{whitt2002} for the definition
(page 81) of the $M_1$ distance and for the properties of the $M_1$ and $J_1$
Skorohod topologies we use below.

Recall that for all $s>0$, $M_s$ denote the number of complete cycles initiated
after time~$0$, and finishing before time $s$. To simplify the
notation, we
assume that $\esp[\phi(X_h(0))]=0$, i.e. that $\phi=\bar\phi$.
Define the
following processes:
\begin{eqnarray*}
\mathcal S_T(u) &=& \frac1{a(T)}\sum_{j=1}^{[Tu]} Z_j( \phi)  , \qquad
 \xi_T(u) = \frac1{a(T)}(M_{Tu}-Tu/\esp[C_1])   ,
\\
\tilde{\mathcal S}_T(u) &=& \mathcal{S}_T(M_{Tu}/T) = \frac1{a(T)}
\sum_{j=1}^{M_{Tu}} Z_j( \phi)
 , \\
R_{0,T} &=& \frac1{a(T)} \int_0^{S_0} \phi(X_h(s)) \, \rmd s  , \qquad
   R_T(u) = \frac1{a(T)} \int_{S_{M_{Tu}}}^{Tu} \phi(X_h(s))
\, \rmd s   .
\end{eqnarray*}
Remark that, if $u<S_0$, then $M_u= 0$ and, hence, $\tilde{\mathcal
S}_T(u)=0$ with the convention $\sum_{j=1}^0(\cdots )=0$. Then
\begin{eqnarray*}
\mathcal Z_T(\phi,u) = R_{0,T} + \tilde{\mathcal S}_T(u) + R_T(u)   .
\end{eqnarray*}

We proceed through a sequence of steps. Specifically, we will prove
that, as $T\to\infty$,
\renewcommand\theenumi{(\roman{enumi})}
\renewcommand\labelenumi{\theenumi}
\begin{enumerate}[(iii)]
\item\label{itemconvergence-sommes-discretes} $S_T$ converges
weakly in $\mathcal D([0,\infty))$ endowed with
the $J_1$ topology to the L\'evy $\alpha$-stable process
$(\esp[C_1])^{1/\alpha} \Lambda(\phi,\cdot)$, where $\Lambda$ is
defined by~\eqref{eqlimit-process};
\item\label{itemverwaat-whitt} $\xi_T$ converges weakly in
$\mathcal
D([0,\infty))$ endowed with the $M_1$
topology to an $\alpha$-stable L\'evy process;
\item\label{itemanscombe} $\tilde{\mathcal S}_T$ converges weakly
in $\mathcal D([0,\infty))$ endowed with
the $J_1$ topology to the L\'evy $\alpha$-stable process
$\Lambda(\phi,\cdot)$;
%
\item\label{iteminterpolatioin-M1} $d_{M_1}(\tilde{\mathcal
S}_T,\mathcal Z_T) \to0$ in probability.
\end{enumerate}

The statement of the theorem will follow from statements
\ref{itemanscombe} and \ref{iteminterpolatioin-M1}. It is
interesting that the statement \ref{iteminterpolatioin-M1}
holds even though $R_T$ converges to zero in neither of the Skorohod
topologies, since otherwise it would then converge uniformly
(because convergence in one of these topology to a continuous limit
implies uniform convergence), and this would imply that $\mcz_T$
weakly converges in the $J_1$ topology to its limit, which is not
possible since the limit is not continuous.

We now prove~\ref{itemconvergence-sommes-discretes}. In the case
$h=0$, the
random variables $Z_j(\phi)$ are i.i.d., centered and their tail
behavior is
given by Lemma~\ref{lemjoint-tails}. The weak convergence in the space
$\mathcal D$ endowed with the $J_1$ topology of the normalized partial sum
process $\mathcal S_T$ to the $\alpha$-stable L\'evy process
$(\esp[C_1])^{1/\alpha}\Lambda(\phi,\cdot)$ is well known in this
case; see,
for example, Resnick \cite{resnick2007}, Corollary~7.1. When $h>0$, $\{
Z_j(\phi)\}$ is no
longer an i.i.d. sequence, so we use the following decomposition. For
$j\geq1$, we write $Z_j(\phi) = Z_{1,j} + Z_{2,j}$ with
\[
Z_{1,j} = \int_{S_{j-1}}^{(S_j-h)\vee S_{j-1}} \phi(X_h(s))
\,\rmd s - \esp \biggl[\int_{S_{j-1}}^{(S_j-h)\vee S_{j-1}}
\phi(X_h(s)) \, \rmd s  \biggr]  .
\]
Observe that the sequence $\{Z_{1,j}\}$ is i.i.d. and centered,
while the sequence
$\{Z_{2,j}\}$ is centered and exponentially $\alpha$-mixing by
(\ref{eqmixing}). Furthermore,
$|Z_{2,j}|\leq 2\|\phi\|_\infty h$. Therefore, by the maximal
inequality for mixing sequences Rio \cite{rio2000}, Theorem~3.1, we
obtain
\[
\esp \Biggl[ \max_{1\leq k\leq n}  \Biggl| \frac1{a(n)} \sum_{j=1}^k
Z_{2,j}  \Biggr|^2  \Biggr] = \mathrm{O}(na_n^{-2}) = \mathrm{o}(1)   .
\]
This implies that the family of processes $a(n)^{-1}\sum
_{j=1}^{[n\cdot]}
Z_{2,j}$ converges weakly to 0 uniformly on compact sets. Since the random
variables $Z_{2,j}$ are uniformly bounded, $Z_{1,j}$ has the same
tail behaviour as
$Z_j$. Thus, as in the case $h=0$, the family of processes
$a(n)^{-1}\sum_{j=1}^{[n\cdot]} Z_{1,j}$ converges weakly in the space
$\mathcal D$ endowed with the $J_1$ topology to the $\alpha$-stable L\'evy
process $(\esp[C_1])^{1/\alpha}\Lambda(\phi,\cdot)$. This
proves~\ref{itemconvergence-sommes-discretes}.

By the regenerative property of the cycles and Lemma~\ref{lemtailC}, $M_t$
is the counting process associated with a renewal process whose interarrival
times $C_j$ are in the domain of attraction of a stable law with index
$\alpha$. More specifically, by Lemmas \ref{lemtailC} and
\ref{lemjoint-tails}, the tails of $C_1$ and $Z_1(\phi)$ are
equivalent. Now~\ref{itemverwaat-whitt} follows from Whitt \cite{whitt2002}, Theorem~4.5.3 and Theorem~6.3.1.

We now prove~\ref{itemanscombe} by the $J_1$-continuity of composition
argument. Observe that $\tilde{\mathcal S}_T=\mathcal S_T\circ
[M_{T\cdot}/T]$. Moreover, $M_{T u}/T=a(T)\xi_T(u)/T+u/\esp[C_1]$
for all
$u\geq0$. Since the supremum functional is continuous in the $M_1$ topology
and $a(T)/T\to0$, we can use \ref{itemverwaat-whitt} to see that $M_{T
\cdot}/T$ converges in the uniform topology on compact intervals to the
linear function $\cdot/\esp[C_1]$ in
probability. By \ref{itemconvergence-sommes-discretes} and Theorem
4.4 in
Billingsley \cite{billingsley1968} we conclude that $(\mathcal S_T ,M_{T\cdot}/T)$
converges weakly to
$ ((\esp[C_1])^{1/\alpha}\Lambda(\phi,\cdot),\cdot/\esp
[C_1] )$ in
the product space $\mathcal D([0,\infty))\times\mathcal D([0,\infty
))$, where
each of the components is endowed with the~$J_1$ topology on compact
intervals. Since the linear function is continuous and strictly increasing,
we can use Theorem~13.2.2 in Whitt \cite{whitt2002} to conclude that
$\tilde{\mathcal S}_T$ converges weakly to
$ (\esp[C_1])^{1/\alpha}\Lambda(\phi,\cdot/\esp[C_1] )$ in
$D([0,\infty))$ endowed with the $J_1$ topology. By the
self-similarity of
centered L\'evy stable motions, the latter process has the same law as~$\Lambda(\phi,\cdot)$. This gives~\ref{itemanscombe}.

It remains to prove \ref{iteminterpolatioin-M1}. Define the process
$\tilde
{\mathcal Z}_T$ by
\[
\tilde{\mathcal Z}_T(t) = \mathcal Z_T(\phi, t) - \mathcal Z_T
(\phi, S_0/T) = a(T)^{-1}\int_{S_{0}}^{Tt} \phi(X_h(s)) \, \rmd s   .
\]
Then, since $S_0<\infty$ a.s.,
\begin{eqnarray*}
\|\tilde{\mathcal Z}_T - \mathcal Z_T\|_\infty=  \biggl|\frac1{a(T)}
\int_0^{S_0} \phi(X_h(s))  \biggr| \leq
\frac{\|\phi\|_\infty S_0}{a(T)} = \mathrm{o}_P(1)  .
\end{eqnarray*}
Since $\tilde{\mathcal S}_T(t)=0$ for all $t\in[0,S_0/T]$, we also have
\[
\sup_{t\in[0,S_0/T]} |\tilde{\mathcal Z}_T(t) - \tilde
{\mathcal
S}_T(t) |\leq\frac{\|\phi\|_\infty S_0}{a(T)}  .
\]
Next, we partition the random interval $[0,S_{M_{T}+1}/T]\supseteq
[0,1]$ into the adjacent intervals
\[
[0,S_0/T]\cup[S_{0}/T,S_1/T]\cup\cdots \cup[S_{i-1}/T,S_i/T]\cup\cdots \cup[S_{M_{T}}/T,S_{M_{T}+1}/T]   .
\]
Recall the following property of the $M_1$ metric: if
$a<b<c$ and
$x_1,x_2$ are functions in $D ([a,c] )$, then
\[
d_{M_1} ( x_1,x_2,[a,c] )\leq\max [ d_{M_1} (
x_1,x_2,[a,b] ),  d_{M_1} ( x_1,x_2,[b,c] ) ] .
\]
We conclude that
\begin{eqnarray*}
d_{M_1}(\tilde{\mathcal S}_T,\mathcal Z_T) & \leq&
d_{M_1}(\mathcal Z_T,\tilde{\mathcal Z}_T) +
d_{M_1}(\tilde{\mathcal Z}_T, \tilde{\mathcal S}_T) \\
& \leq&\frac{2\|\phi\|_\infty S_0}{a(T)} + \max_{i=1,\ldots,M_{T}}
d_{M_1}(\tilde{\mathcal
Z}_T , \tilde{\mathcal S}_T, [S_{i-1}/T,S_i/T]) \\
&&{} + d_{M_1}(\tilde{\mathcal
Z}_T , \tilde{\mathcal S}_T, [S_{M_T}/T,1])  .
\end{eqnarray*}
Notice that the last term in the right-hand side is bounded by $\|\phi
\|_\infty
C_{M_T+1}/a(T)$, and the finite mean of $C_1$ implies that the $C_{M_T+1}$
converges weakly as $T\to\infty$ and, in particular, the family of
the laws of
$(C_{M_T+1})$ is tight. Observe, further, that $\tilde{\mathcal Z}_T$
continuously interpolates $\tilde{\mathcal S}_T$ at the points $t=S_i/T$,
$i=0,1,2,\ldots.$ Hence, by  (\ref{eqMTbound}), $\pr(T>S_0)\to1$ and
stationarity we see that for any $\eta>0$,
\[
\pr \bigl(d_{M_1}(\tilde{\mathcal S}_T,\mathcal Z_T) >
\eta \bigr) \leq\frac{2T}{\esp[C_1]} \pr \bigl(
d_{M_1}(\tilde{\mathcal Z}_T, \tilde{\mathcal S}_T, [S_{0}/T,S_1/T])
> \eta/2  \bigr) + \mathrm{o}(1)   .
\]
Henceforth, we now only consider the process $X_h(t)$ on $[S_0, S_1]$.
We use
the notation introduced in Section~\ref{secnotat-prel-results}. First
of all,
\begin{eqnarray*}
d_{M_1}(\tilde{\mathcal Z}_T,\tilde{\mathcal
S}_T,[S_{0}/T,S_1/T])&\leq&
\sup_{u\in[S_{0}/T,S_1/T]} |\tilde{\mathcal
Z}_T(u)-\tilde{\mathcal S}_T(u) |
\\[-3pt]
&\leq&
a(T)^{-1}\sup_{v\in[S_0,S_1]}\int_{S_{0}}^{v} \phi(X_h(s)) \,
\rmd s \leq a(T)^{-1}C_1\|\phi\|_\infty .
\end{eqnarray*}
Combining this with (\ref{eqC1Nnul}), we see that for any $\eta>0$,
\[
\pr \bigl(d_{M_1}(\tilde{\mathcal Z}_T,\tilde{\mathcal
S}_T,[S_{0}/T,S_1/T])>\eta; N_{\epsilon,T}=0 \bigr)=\mathrm{o}(T^{-1}) ,
\]
as long as $\epsilon>0$ is chosen to be small enough.

Next, we consider the event $\{N_{\epsilon,T}\geq1\}$.
Define
\[
\check{\mathcal Z}_T(t) = a(T)^{-1} \int_{S_0}^{t T}
\mce(W_{\epsilon,T}, \phi) \1_{[\Gamma_{\epsilon,T},
\Gamma_{\epsilon,T} + Y_{\epsilon,T})}(s) \, \rmd s .
\]
Observe that $\check{\mathcal Z}_T$ is monotone on $[S_0/T,S_1/T]$ and
piecewise
linear and $\tilde{\mathcal S}_T$ is constant on $[S_0/T,S_1/T)$ with
a step at
the point $S_1/T$. Using these properties and the definition of the $M_1$
distance, it is not difficult to check that
\[
d_{M_1}(\check{\mathcal Z}_T, \tilde{\mathcal S}_T,[S_{0}/T,S_1/T])
\leq\frac{C_1}{T} \vee | \tilde{\mathcal S}_T(S_1/T) -
\check{\mathcal Z}_T(S_1/T)  |   .
\]
On the other hand, bounding by the uniform distance gives us
\[
d_{M_1}(\tilde{\mathcal Z}_T,\check{\mathcal
Z}_T,[S_{0}/T,S_1/T])\leq\sup_{t\in[S_{0}/T,S_1/T]}|\tilde{\mathcal
Z}_T(t)-\check{\mathcal Z}_T(t)| .
\]
Since $\tilde{\mathcal S}_T(S_1/T)=\tilde{\mathcal Z}_T(S_1/T)$, the previous
bounds yield
\begin{eqnarray*}
&&\pr\bigl( d_{M_1}(\tilde{\mathcal Z}_T, \tilde{\mathcal
S}_T,[S_{0}/T,S_1/T]) >
\eta; N_{\epsilon,T}=1 \bigr) \\[-2pt]
&& \quad \leq\pr( C_1 > \eta T/2; N_{\epsilon,T}=1 ) + 2 \pr\Bigl(
\sup_{t\in[S_{0}/T,S_1/T]}|\tilde{\mathcal Z}_T(t)-\check{\mathcal
Z}_T(t)| >
\eta/2; N_{\epsilon,T}=1 \Bigr).
\end{eqnarray*}
By Lemma~\ref{lemtailC}, we know that $\pr(C_1>\eta T)=\mathrm{o}(T^{-1})$.
Moreover, since
\begin{eqnarray*}
&&\sup_{t\in[S_{0}/T,S_1/T]} | \tilde{\mathcal Z}_T(t) -
\check{\mathcal Z}_T(t)|\\[-2pt]
 && \quad = \frac1{a(T)} \sup_{v\in[S_0,S_1]}  \biggl|
\int_{S_{0}}^{v} \bigl\{ \phi(X_h(s)) - \mce(W_{\epsilon,T},
\phi) \1_{[\Gamma_{\epsilon,T},\Gamma_{\epsilon,T} +
Y_{\epsilon,T})}(s)\bigr\} \, \rmd s \biggr|   ,
\end{eqnarray*}
Lemma~\ref{lemM1} states exactly that
\begin{eqnarray*}
\pr \Bigl( \sup_{t\in[S_{0}/T,S_1/T]}|\tilde{\mathcal Z}_T(t) -
\check{\mathcal
Z}_T(t)| > \eta; N_{\epsilon,T}\geq1 \Bigr) = \mathrm{o}(T^{-1})   .
\end{eqnarray*}
This completes the proof
of~\ref{iteminterpolatioin-M1}.\vadjust{\goodbreak}
\end{pf*}

\section*{Acknowledgements}
The authors would like to thank the two
referees and an associate editor for their unusually
careful reading of the paper and helpful anonymous comments.
Roueff and Soulier's research were partially supported by
the ANR Grant ANR-09-BLAN-0029-01. Samorodnitsky's research was
partially supported by the ARO
Grant W911NF-07-1-0078 at Cornell University, Department of
Mathematics, Universit\'e de Paris Ouest Nanterre during his visit in
2007, Laboratory of
Actuarial Mathematics, Department of Mathematics, University of
Copenhagen and by Laboratory of Informatics, Technical University of
Denmark, during his sabbatical stay in 2008--2009.


\printhistory


\begin{thebibliography}{20}

\bibitem{billingsley1968}
\begin{bbook}[mr]
\bauthor{\bsnm{Billingsley},~\bfnm{Patrick}\binits{P.}}
(\byear{1968}).
\btitle{Convergence of Probability Measures}.
\baddress{New York}: \bpublisher{Wiley}.
\bid{mr={0233396}}
\end{bbook}
\endbibitem

\bibitem{cohn1972}
\begin{barticle}[mr]
\bauthor{\bsnm{Cohn},~\bfnm{Donald~L.}\binits{D.L.}}
(\byear{1972}).
\btitle{Measurable choice of limit points and the existence of separable and
  measurable processes}.
\bjournal{Z. Wahrsch. Verw. Gebiete}
\bvolume{22}
\bpages{161--165}.
\bid{mr={0305444}}
\end{barticle}
\endbibitem

\bibitem{dehlingtaqqu1989}
\begin{barticle}[mr]
\bauthor{\bsnm{Dehling},~\bfnm{Herold}\binits{H.}} \AND
  \bauthor{\bsnm{Taqqu},~\bfnm{Murad~S.}\binits{M.S.}}
(\byear{1989}).
\btitle{The empirical process of some long-range dependent sequences with an
  application to {$U$}-statistics}.
\bjournal{Ann. Statist.}
\bvolume{17}
\bpages{1767--1783}.
\bid{doi={10.1214/aos/1176347394}, issn={0090-5364}, mr={1026312}}
\end{barticle}
\endbibitem

\bibitem{fayroueffsoulier2007}
\begin{barticle}[mr]
\bauthor{\bsnm{Fay},~\bfnm{Gilles}\binits{G.}},
  \bauthor{\bsnm{Roueff},~\bfnm{Fran{\c{c}}ois}\binits{F.}} \AND
  \bauthor{\bsnm{Soulier},~\bfnm{Philippe}\binits{P.}}
(\byear{2007}).
\btitle{Estimation of the memory parameter of the infinite-source {P}oisson
  process}.
\bjournal{Bernoulli}
\bvolume{13}
\bpages{473--491}.
\bid{doi={10.3150/07-BEJ5123}, issn={1350-7265}, mr={2331260}}
\end{barticle}
\endbibitem

\bibitem{hall1988}
\begin{bbook}[mr]
\bauthor{\bsnm{Hall},~\bfnm{Peter}\binits{P.}}
(\byear{1988}).
\btitle{Introduction to the Theory of Coverage Processes}.
\bseries{Wiley Series in Probability and Mathematical Statistics: Probability
  and Mathematical Statistics}.
\baddress{New York}: \bpublisher{Wiley}.
\bid{mr={0973404}}
\end{bbook}
\endbibitem

\bibitem{heffernanresnick2007}
\begin{barticle}[mr]
\bauthor{\bsnm{Heffernan},~\bfnm{Janet~E.}\binits{J.E.}} \AND
  \bauthor{\bsnm{Resnick},~\bfnm{Sidney~I.}\binits{S.I.}}
(\byear{2007}).
\btitle{Limit laws for random vectors with an extreme component}.
\bjournal{Ann. Appl. Probab.}
\bvolume{17}
\bpages{537--571}.
\bid{doi={10.1214/105051606000000835}, issn={1050-5164}, mr={2308335}}
\end{barticle}
\endbibitem

\bibitem{kallenberg2002}
\begin{bbook}[mr]
\bauthor{\bsnm{Kallenberg},~\bfnm{Olav}\binits{O.}}
(\byear{2002}).
\btitle{Foundations of Modern Probability},
\bedition{2nd} ed.
\bseries{Probability and Its Applications (New York)}.
\baddress{New York}: \bpublisher{Springer}.
\bid{mr={1876169}}
\end{bbook}
\endbibitem

\bibitem{maulikresnickrootzen2002}
\begin{barticle}[mr]
\bauthor{\bsnm{Maulik},~\bfnm{Krishanu}\binits{K.}},
  \bauthor{\bsnm{Resnick},~\bfnm{Sidney}\binits{S.}} \AND
  \bauthor{\bsnm{Rootz{\'e}n},~\bfnm{Holger}\binits{H.}}
(\byear{2002}).
\btitle{Asymptotic independence and a network traffic model}.
\bjournal{J. Appl. Probab.}
\bvolume{39}
\bpages{671--699}.
\bid{issn={0021-9002}, mr={1938164}}
\end{barticle}
\endbibitem

\bibitem{mikoschresnickrootzenstegeman2002}
\begin{barticle}[mr]
\bauthor{\bsnm{Mikosch},~\bfnm{Thomas}\binits{T.}},
  \bauthor{\bsnm{Resnick},~\bfnm{Sidney}\binits{S.}},
  \bauthor{\bsnm{Rootz{\'e}n},~\bfnm{Holger}\binits{H.}} \AND
  \bauthor{\bsnm{Stegeman},~\bfnm{Alwin}\binits{A.}}
(\byear{2002}).
\btitle{Is network traffic approximated by stable {L}\'evy motion or fractional
  {B}rownian motion?}
\bjournal{Ann. Appl. Probab.}
\bvolume{12}
\bpages{23--68}.
\bid{doi={10.1214/aoap/1015961155}, issn={1050-5164}, mr={1890056}}
\end{barticle}
\endbibitem

\bibitem{resnick1992}
\begin{bbook}[mr]
\bauthor{\bsnm{Resnick},~\bfnm{Sidney I.}\binits{S.I.}}
(\byear{1992}).
\btitle{Adventures in Stochastic Processes}.
\baddress{Boston, MA}: \bpublisher{Birkh\"auser}.
\bid{mr={1181423}}
\end{bbook}
\endbibitem

\bibitem{resnick2007}
\begin{bbook}[mr]
\bauthor{\bsnm{Resnick},~\bfnm{Sidney~I.}\binits{S.I.}}
(\byear{2007}).
\btitle{Heavy-Tail Phenomena: Probabilistic and Statistical Modeling}.
\bseries{Springer Series in Operations Research and Financial Engineering}.
\baddress{New York}: \bpublisher{Springer}.
\bid{mr={2271424}}
\end{bbook}
\endbibitem

\bibitem{resnicksamorodnitsky1999a}
\begin{barticle}[mr]
\bauthor{\bsnm{Resnick},~\bfnm{Sidney I.}\binits{S.I.}} \AND
  \bauthor{\bsnm{Samorodnitsky},~\bfnm{Gennady}\binits{G.}}
(\byear{1999}).
\btitle{Activity periods of an infinite server queue and performance of certain
  heavy tailed fluid queues}.
\bjournal{Queueing Systems Theory Appl.}
\bvolume{33}
\bpages{43--71}.
\bid{issn={0257-0130}, mr={1748638}}
\end{barticle}
\endbibitem

\bibitem{resnickvandenberg2000a}
\begin{barticle}[mr]
\bauthor{\bsnm{Resnick},~\bfnm{Sidney I.}\binits{S.I.}} \AND
  \bauthor{\bparticle{van~den} \bsnm{Berg},~\bfnm{Eric}\binits{E.}}
(\byear{2000}).
\btitle{Weak convergence of high-speed network traffic models}.
\bjournal{J.~Appl. Probab.}
\bvolume{37}
\bpages{575--597}.
\bid{issn={0021-9002}, mr={1781014}}
\end{barticle}
\endbibitem



\bibitem{rio2000}
\begin{bbook}[mr]
\bauthor{\bsnm{Rio},~\bfnm{Emmanuel}\binits{E.}}
(\byear{2000}).
\btitle{Th\'eorie Asymptotique des Processus Al\'eatoires Faiblement
  D\'ependants}.
\bseries{Math\'ematiques \& Applications (Berlin) [Mathematics \&
  Applications]}
\bvolume{31}.
\baddress{Berlin}: \bpublisher{Springer}.
\bid{mr={2117923}}
\end{bbook}
\endbibitem

\bibitem{samorodnitskytaqqu1994}
\begin{bbook}[mr]
\bauthor{\bsnm{Samorodnitsky},~\bfnm{Gennady}\binits{G.}} \AND
  \bauthor{\bsnm{Taqqu},~\bfnm{Murad~S.}\binits{M.S.}}
(\byear{1994}).
\btitle{Stable Non-{G}aussian Random Processes: Stochastic Models with Infinite Variance}.
\bseries{Stochastic Modeling}.
\baddress{New York}: \bpublisher{Chapman \& Hall}.
\bid{mr={1280932}}
\end{bbook}
\endbibitem

\bibitem{surgailis2002}
\begin{barticle}[mr]
\bauthor{\bsnm{Surgailis},~\bfnm{Donatas}\binits{D.}}
(\byear{2002}).
\btitle{Stable limits of empirical processes of moving averages with infinite
  variance}.
\bjournal{Stochastic Process. Appl.}
\bvolume{100}
\bpages{255--274}.
\bid{doi={10.1016/S0304-4149(02)00103-5}, issn={0304-4149}, mr={1919616}}
\end{barticle}
\endbibitem

\bibitem{surgailis2004}
\begin{barticle}[mr]
\bauthor{\bsnm{Surgailis},~\bfnm{Donatas}\binits{D.}}
(\byear{2004}).
\btitle{Stable limits of sums of bounded functions of long-memory moving
  averages with finite variance}.
\bjournal{Bernoulli}
\bvolume{10}
\bpages{327--355}.
\bid{doi={10.3150/bj/1082380222}, issn={1350-7265}, mr={2046777}}
\end{barticle}
\endbibitem

\bibitem{vonbahresseen1965}
\begin{barticle}[mr]
\bauthor{\bparticle{von} \bsnm{Bahr},~\bfnm{Bengt}\binits{B.}} \AND
  \bauthor{\bsnm{Esseen},~\bfnm{Carl-Gustav}\binits{C.G.}}
(\byear{1965}).
\btitle{Inequalities for the {$r$}th absolute moment of a~sum of random
  variables, {$1\leq r\leq 2$}}.
\bjournal{Ann. Math. Statist}
\bvolume{36}
\bpages{299--303}.
\bid{issn={0003-4851}, mr={0170407}}
\end{barticle}
\endbibitem

\bibitem{whitt2002}
\begin{bbook}[mr]
\bauthor{\bsnm{Whitt},~\bfnm{Ward}\binits{W.}}
(\byear{2002}).
\btitle{Stochastic-Process Limits: An Introduction to Stochastic-Process Limits and Their Application to
  Queues}.
\bseries{Springer Series in Operations Research}.
\baddress{New York}: \bpublisher{Springer}.
\bid{mr={1876437}}
\end{bbook}
\endbibitem

\end{thebibliography}
\end{document}